\documentclass[a4paper,leqno,11pt]{article}
\usepackage{amssymb,amsfonts,amsmath,amsthm}
\usepackage{epsfig,epsf,curves}
\usepackage{bbm,epic,eepic}
\usepackage[mathscr]{eucal}
\theoremstyle{remark}

\begin{document}

\sloppy


\def\Mat#1#2#3#4{\left(\!\!\!\begin{array}{cc}{#1}&{#2}\\{#3}&{#4}\\ \end{array}\!\!\!\right)}

\def\tAA{{\Bbb A}}
\def\CC{{\Bbb C}}
\def\HH{{\Bbb H}}
\def\NN{{\Bbb N}}
\def\QQ{{\Bbb Q}}
\def\RR{{\Bbb R}}
\def\ZZ{{\Bbb Z}}
\def\FF{{\Bbb F}}
\def\SS{{\Bbb S}}
\def\GG{{\Bbb G}}
\def\PP{{\Bbb P}}
\def\LL{{\Bbb L}}


\title{On truncation of irreducible representations of Chevalley groups II}
\author {J. Mahnkopf}
\maketitle

{\small {\bf Abstract} We determine the growth of the dimension of the slope subspaces of the cohomology of arithmetic subgroups 
in reductive algebraic groups as a function of the slope.

}

\bigskip

{\bf Introduction. } {\bf (0.1)} 
%
%
According to part of (a higher rank analogue of) the Gouvea-Mazur Conjecture, the dimension of the slope $\alpha$-subspace of the 
cohomology of an arithmetic subgroup $\Gamma$ in a reductive group $\tilde{\bf G}$ with coefficients in a finite dimensional 
irreducible representation $L$ is bounded independently of $L$. In the original case ${\bf GL}_2$ this has been proven by 
Wan (following work of Coleman) and - using different methods - by Buzzard (cf. [B], [W]). Moreover, these authors determine explicit bounds, 
i.e. they give an explicit linear polynomial in $\alpha$ which yields an upper bound for the dimension of the slope 
$\alpha$ subspace of the cohomology (cf. [B], Corollary 1 in sec. 3).

\medskip

In this article we determine bounds for the dimension of the slope subspaces of the cohomology of arithmetic subgroups in the higher rank case. 
To state our result, we fix a prime $p\in{\Bbb N}$. We let $\tilde{\bf G}$ be a ${\Bbb Q}$-split reductive group with maximal ${\Bbb Q}$-split 
torus $\tilde{\bf T}$, $\Gamma\le \tilde{\bf G}({\Bbb Q})$ an arithmetic subgroup satisfying a certain level condition at $p$ and we denote 
by $L_{\tilde{\lambda}}$ the irreducible $\tilde{\bf G}$-module of highest weight $\tilde{\lambda}\in X(\tilde{\bf T})$. Moreover, $s=\#\,\Phi^+$ 
is the number of positive roots of $\tilde{\bf G}$. Finally, $H^i(\Gamma,L_{\tilde{\lambda}})^{\le\alpha}$ denotes the slope $\le\alpha$-subspace of 
the cohomology of $\Gamma$ w.r.t. a normalized Hecke operator at $p$.

\bigskip

{\bf Theorem } (cf. (2.3) Corollary). {\it There are natural numbers $m,n$ only depending on $\Gamma$ such that 
$$
{\rm dim}\,H^i(\Gamma,L_{\tilde{\lambda}})^{\le\alpha}\le m\alpha^s+n
$$
for all $\alpha\in{\Bbb Q}_{\ge 0}$, all dominant and algebraic weights $\tilde{\lambda}$ and all $i$.
}

\bigskip

Thus, in the case ${\bf GL}_2$, as in [B], we obtain a linear polynomial in $\alpha$. 
On the other hand, in [B] the factor $m$ is made explicit by showing that it can be chosen equal to the cardinality of a minimal system 
of generators for $\Gamma$. 
%
%
In the higher rank case, we do not know whether 
the coefficients $m$ and $n$ have an interpretation in terms of $\Gamma$. We show that they can be computed explicitly from 
the maximum of the numbers $n_i$, where $n_i$ is the number of $i$-cells of a cell decomposition ${\cal Z}$ of the locally symmetric space $\Gamma\backslash \bar{X}$ and from the coefficients of the Bernoulli polynomials $B_s$ and $B_{s+1}$ (cf. (1.4) and (2.2)). 
By Tietze's Theorem, the $1$-cells in ${\cal Z}$ yield a finite system of generators for $\Gamma$, hence, the determination of the number of $1$ cells (which is sufficient in the ${\bf GL}_2$-case) and of a finite 
system of generators for $\Gamma$ are related problems. In this sense, the maximum of the $n_i$'s may be seen as a higher rank analogue of 
the number of generators of $\Gamma$.

\bigskip

{\bf (0.2) } The proof of the Theorem relies on the methods in [M], where the existence of a bound for the dimension of the slope subspaces which is independent 
of the coefficient system $L$ has been proven. In addition, we use a Theorem of Borel-Serre according 
to which any arithmetic group is of type FL. This also simplifies the arguments in [M] since the result
about isomorphisms between truncations of representations is no longer necessary. On the other hand, congruences between truncations might 
also yield local constancy of the dimension of the slope subspaces (cf. [P] in the case of unit groups of quaternion algebras) which may not be 
possible using the simpler argument. Following [B], in order to obtain bounds for the dimension 
of the slope subspaces, we determine lower bounds for the Newton polygon of the normalized Hecke operator at $p$ acting on cohomology.

\medskip

We hope that the Theorem will have applications to the construction of explicit congruences between automorphic eigenforms.



\bigskip

{\bf (0.3) Notations. } We keep the notations and definitions from [M]. For convenience, we recall the most important ones, referring 
for the remaining ones to [M]. $\tilde{\bf G}$ is a reductive, ${\Bbb Q}$-split algebraic group with derived group ${\bf G}$ and 
$\Gamma\le{\bf G}({\Bbb Q})$ is an arithmetic subgroup. We fix a prime $p\in{\Bbb N}$ and we assume that $\Gamma$ is contained in the 
level subgroup $K_*(p)\le{\bf G}({\Bbb Z}_p)$ defined in [M], sec. 4.2 ($K_*(p)$ is a Iwahori type subgroup). We denote by $\tilde{\bf T}$ 
resp. ${\bf T}$ a maximal split torus in $\tilde{\bf G}$ resp. in ${\bf G}$ such that ${\bf T}\le\tilde{\bf T}$. We write 
${\mathfrak g}$ resp. ${\mathfrak h}$ for the Lie Algebra of ${\bf G}$ resp. ${\bf T}$ and for any 
weight $\tilde{\lambda}\in X(\tilde{\bf T})$ we set $\lambda={\rm d}\,\tilde{\lambda}|_{\bf T}\in{\mathfrak h}^*$. Let 
$\tilde{\lambda}\in X(\tilde{\bf T})$ be an algebraic and dominant weight. We denote by $L_{\tilde{\lambda}}$ the irreducible 
representation of $\tilde{\bf G}$ of highest weight $\tilde{\lambda}$ and by $V_\lambda={\cal U}^-{\rm v}_\lambda$ resp. 
$L_\lambda={\cal U}^-v_\lambda$ the Verma module resp. the irreducible ${\mathfrak g}$-module of highest weight 
$\lambda={\rm d}\,\tilde{\lambda}|_{\bf T}$. Thus, $L_\lambda$ is the derived representation of the restriction of $L_{\tilde{\lambda}}$ to ${\bf G}$. 
Moreover, $\Pi_\lambda\subseteq{\mathfrak h}^*$ is the set of weights occuring in $L_\lambda$ and $\Phi^+$ is the set of positive roots of 
${\mathfrak g}$. 

\medskip

Finally, we define the $\Gamma$-modules $L_\lambda({\Bbb Z}_p,r)$, $L_{\tilde{\lambda}}({\Bbb Z}_p,r)$ and 
${\bf L}^{[r]}_\lambda({\Bbb Z}_p)=L_\lambda({\Bbb Z}_p)/L_\lambda({\Bbb Z}_p,r)$ as in [M]. The modules $L_\lambda({\Bbb Z}_p,r)$ and 
$L_{\tilde{\lambda}}({\Bbb Z}_p,r)$ are isomorphic but on $L_{\tilde{\lambda}}({\Bbb Z}_p,r)$ (like on $L_{\tilde{\lambda}}({\Bbb Z}_p)$) we still 
have an action of dominant elements $h\in\tilde{\bf T}({\Bbb Q}_p)$ (cf. [M], 5.4 Lemma). We fix a strictly dominant element 
$h\in\tilde{\bf T}({\Bbb Q})$; the normalized Hecke operator ${\Bbb T}(h)$ (cf. [M], sec. 5.3) then acts on cohomology with 
integral coefficients $H^i(\Gamma,L_{\tilde{\lambda}}({\Bbb Z}_p))$, $H^i(\Gamma,L_{\tilde{\lambda}}({\Bbb Z}_p,r))$ and 
$H^i(\Gamma,{\bf L}^{[r]}_{\tilde{\lambda}}({\Bbb Z}_p))$ (cf. [M], 5.8 Proposition).

\vspace{1cm}

\section{The Newton polygon}


\bigskip

{\bf (1.1) } We recall that we have chosen an ordering of the positive roots $\Phi^+=\{\alpha_1,\ldots,\alpha_s\}$ of ${\bf G}$. We use 
this to define the following counting function. We denote by $N_h$, $h\in{\Bbb N}_0$, the number of tuples
${\bf n}=(n_1,\ldots,n_s)\in{\Bbb N}_0^s$ such that
$$
\sum_{i=1}^s n_i{\rm ht}(\alpha_i)=h.
$$
The above condition implies that any coordinate $n_i$ of ${\bf n}$ is smaller than or equal to $h$ 
and that the last coordinate $n_s$ is determined by the first $s-1$ coordinates. Hence, we obtain $N_h\le (h+1)^{s-1}$. 

\bigskip

{\bf (1.2) Lemma. } {\it For any dominant and algebraic weight $\lambda$ and any $r\in{\Bbb N}$ there is an embedding
$$
{\bf L}_\lambda^{[r]}({\Bbb Z}_p)\le\bigoplus_{h=0}^r ({\Bbb Z}_p/{p^{r-h}{\Bbb Z}_p})^{N_h}.
$$
}

\medskip

{\it Proof. } By equation (5) in section (2.1) of [M] we know  
$$
V_\lambda({\Bbb Z},\mu)=\bigoplus_{{\bf n}\atop n_1\alpha_1+\cdots+n_s\alpha_s=\lambda-\mu} {\Bbb Z} X_-^{\bf n} {\rm v}_\lambda
$$
Since further $X_-^{\bf n}{\rm v}_\lambda$ has weight $\mu=\lambda-\sum_i n_i\alpha_i$ which has relative height ${\rm ht}_\lambda(\mu)=\sum_i n_i{\rm ht}(\alpha_i)$ we obtain
\begin{eqnarray*}
\bigoplus_{{\bf n}\in{\Bbb N}_0^s\atop {\rm ht}(n_1\alpha_1+\cdots+n_s\alpha_s)=h}{\Bbb Z}X_-^{\bf n} {\rm v}_\lambda
&=&\bigoplus_{\mu\in\Pi_\lambda\atop {\rm ht}_\lambda(\mu)=h}
\bigoplus_{{\bf n}\in{\Bbb N}_0^s\atop \lambda-n_1\alpha_1+\cdots+n_s\alpha_s=\mu} {\Bbb Z} X_-^{\bf n}{\rm v}_\lambda\\
&=&\bigoplus_{\mu\in\Pi_\lambda\atop {\rm ht}_\lambda(\mu)=h} V_\lambda({\Bbb Z},\mu).\\
\end{eqnarray*}
%
%
In particular, we obtain
$$
N_h={\rm dim}_{\Bbb Z}\, \bigoplus_{\mu\in\Pi_\lambda\atop {\rm ht}_\lambda(\mu)=h} V_\lambda({\Bbb Z},\mu).
$$
Since $L_\lambda$ is a quotient of the Verma module $V_\lambda$ this implies
$$
{\rm dim}\, \bigoplus_{\mu\in\Pi_\lambda\atop {\rm ht}_\lambda(\mu)=h} L_\lambda({\Bbb Z},\mu)\le N_h.
$$
By [M], (1.5) equation (3) we know that
\begin{eqnarray*}
{\bf L}_\lambda^{[r]}({\Bbb Z}_p)=\frac{L_\lambda({\Bbb Z}_p)}{L_\lambda({\Bbb Z}_p,r)}
&\cong&\bigoplus_{h=0}^r\bigoplus_{\mu\in\Pi_\lambda\atop{\rm ht}_\lambda(\mu)=h} 
\frac{L_\lambda({\Bbb Z}_p,\mu)}{p^{r-h} L_\lambda({\Bbb Z}_p,\mu)}\\
&=&\bigoplus_{h=0}^r \bigoplus_{\mu\in\Pi_\lambda\atop{\rm ht}_\lambda(\mu)=h} L_\lambda({\Bbb Z}_p,\mu)\otimes
\frac{{\Bbb Z}_p}{p^{r-h} {\Bbb Z}_p}.\\
\end{eqnarray*}
We obtain
$$
{\bf L}_\lambda^{[r]}({\Bbb Z}_p)\le\bigoplus_{h=0}^r ({\Bbb Z}_p/{p^{r-h}{\Bbb Z}_p})^{N_h}
$$
and the Lemma is proven.

\bigskip

{\bf (1.3) } According to [B-S] the arithmetic subgroup $\Gamma\le {\bf G}({\Bbb Q})$ is a $FL$-group (cf. [B], p. 218). 
Thus, there is a resolution of the trivial $\Gamma$-module ${\Bbb Z}$ 
$$
0\rightarrow M_d\rightarrow \cdots \rightarrow M_1\rightarrow M_0\rightarrow {\Bbb Z}\rightarrow 0
$$
where $M_i$ is a free ${\Bbb Z}\Gamma$-module of finite rank (cf. [B], p. 199). This yields a finite free resolution 
of the trivial $\Gamma$-module ${\Bbb Z}_p$
$$
0\rightarrow M_{p,d}\rightarrow \cdots \rightarrow M_{p,1}\rightarrow M_{p,0}\rightarrow {\Bbb Z}_p\rightarrow 0
$$
where $M_{p,i}=M_i\otimes_{{\Bbb Z}\Gamma}{\Bbb Z}_p\Gamma$; note that ${\Bbb Z}\otimes_{{\Bbb Z}\Gamma}{\Bbb Z}_p\Gamma\cong {\Bbb Z}_p$ as 
${\Bbb Z}_p\Gamma$-modules. 
%
%
%
%
The group $H^\bullet(\Gamma,{\bf L}_\lambda^{[r]}({\Bbb Z}_p))$ then may be computed as the cohomology of the complex
$$
0\rightarrow{\rm Hom}_{{\Bbb Z}_p\Gamma}(M_{p,0},{\bf L}_\lambda^{[r]}({\Bbb Z}_p))
\rightarrow\cdots\rightarrow
{\rm Hom}_{{\Bbb Z}_p\Gamma}(M_{p,d},{\bf L}_\lambda^{[r]}({\Bbb Z}_p))\rightarrow 0;
$$
%
%
%
%
We set
$$
g=g(\Gamma)={\rm sup}_{i=0,\ldots,d} \, {\rm rk}_{{\Bbb Z}\Gamma} M_i={\rm sup}_{i=0,\ldots,d} \, {\rm rk}_{{\Bbb Z}_p\Gamma} M_{p,i}.
$$
Thus, $g$ only depends on the arithmetic group $\Gamma$.

\bigskip

{\bf (1.4) Remark. }  The number $g$ can be related to $\Gamma$ as follows. We denote by $Y=\Gamma\backslash \bar{X}$ the Borel-Serre compactification
of the locally symmetric space $\Gamma\backslash X$. By the work of Borel-Serre $Y$ is a compact $K(\Gamma,1)$ space (cf. [B-S]),
hence, there is a {\it finite} CW complex ${\cal Z}$ which is homtop to $Y$.  
The universal cover
$\tilde{Y}\rightarrow Y$ inherits a structure of $CW$ complex $\tilde{\cal Z}$ from ${\cal Z}$ and the cellular complex $\tilde{C}_\bullet=(\tilde{C}_i)$ attached to $\tilde{\cal Z}$ is a
complex consisting of free ${\Bbb Z}\Gamma$-modules $\tilde{C}_i$. The module $\tilde{C}_i$ has a natural basis which is in bijection with the set of $i$-cells ${\cal Z}$ (cf. [Br], p. 15); hence, $\tilde{C}_\bullet$ is a finite complex consisting of free, finitely generated ${\Bbb Z}\Gamma$-modules.
In particular, the cohomology groups $H^i(\Gamma,{\bf L}_{\lambda}^{[r]}({\Bbb Z}_p))$
can be computed using $\tilde{C}_\bullet$. Since ${\rm rk}_{{\Bbb Z}\Gamma} \,\tilde{C}_i$ equals the number of $i$-cells of ${\cal Z}$ 
this shows that we can take for $g$ the following value: $g={\rm sup}_i n_i$, where $n_i$ is number of $i$-cells of a CW complex ${\cal Z}$ 
which is homotop to $\Gamma\backslash\bar{X}$.


\medskip

If we only look at cohomolgy in degree $1$ (e.g. in the ${\bf GL}_2$-case) then $g$ may be taken as the number of $1$ cells in ${\cal Z}$. Since any 
CW complex for $\Gamma\backslash\bar{X}$ yields a presentation for $\Gamma$ with the generators corresponding to a subset of the set of $1$-cells of 
${\cal Z}$ by Tietze's theorem (cf. [R], Theorem 11.31, p. 374),  
this relates $g$ to the cardinality of finite systems of generators for $\Gamma$.


\bigskip

{\bf (1.5) } We set $t=t(\tilde{\lambda},i)={\rm dim}\, H^i(\Gamma,L_{\tilde{\lambda}}({\Bbb Q}_p))$ and we denote by 
$$
{\cal N}={\cal N}_{\tilde{\lambda},i}:\,[0,t]\rightarrow{\Bbb R}_{\ge 0}
$$ 
the Newton polygon of ${\Bbb T}(h)$ acting on $H^i(\Gamma,L_{\tilde{\lambda}}({\Bbb Q}_p))$. We let $g$ be as in (1.3).


\bigskip

{\bf Theorem. }{\it The Newton polygon ${\cal N}$ lies above the restriction to $[0,t]$ of the 
piecewise linear function $f_\infty:\,{\Bbb R}_{\ge 0}\rightarrow{\Bbb R}_{\ge 0}$ which starts in $(0,0)$ and has slope $0$ in the interval
$$
0\le x\le g\frac{B_s(2)-B_s(0)}{s}
$$
and slope $j$ in the interval
$$
g \frac{B_s(j+1)-B_s(0)}{s}\le x\le g \frac{B_s(j+2)-B_s(0)}{s}, \quad j=1,2,\ldots
$$
Equivalently, $f_\infty$ may be defined as the piecewise linear function joining the points
$$
(0,0),\; P_j=(g \frac{B_s(j+2)-B_s(0)}{s},g \sum_{h=0}^j h (h+1)^{s-1}),\quad j=0,1,2,\ldots.
$$

}

\medskip

{\it Proof. } We proceed in steps.

\medskip

{\bf (1.5.1) } We let $\tilde{\lambda}\in X(\tilde{\bf T})$ and we set $\lambda={\rm d}\,\tilde{\lambda}|_{\bf T}\in{\mathfrak h}^*$. 
Moreover, we select a natural number $r\in{\Bbb N}$. By (1.2) Lemma we know that
$$
{\bf L}_\lambda^{[r]}({\Bbb Z}_p)\le\bigoplus_{h=0}^r ({\Bbb Z}_p/{p^{r-h}{\Bbb Z}_p})^{N_h},
$$
which implies that
$$
{\rm Hom}_{{\Bbb Z}_p\Gamma}(M_{p,i},{\bf L}_\lambda^{[r]}({\Bbb Z}_p))
\le \bigoplus_{h=0}^r  ({\Bbb Z}_p/{p^{r-h}{\Bbb Z}_p})^{N_h {\rm rk}_{{\Bbb Z}_p\Gamma} M_{p,i}}.\leqno(1)
$$
We set $g_i={\rm rk}_{{\Bbb Z}_p\Gamma} M_{p,i}$ and we denote by $(p^{a_l})_l$, $a_1\ge a_2\ge\cdots\ge a_n>0$ the sequence of elementary 
divisors of the right hand side of equation (1), i.e.
$$
(p^{a_l})_{l=1,\ldots,n}=(p^r,\ldots,p^r,p^{r-1},\ldots,p^{r-1},\ldots,p,\ldots,p),\leqno(2)
$$
where $p^{r-h}$ appears $g_i N_h$-times. The exact sequence
$$
0\rightarrow L_\lambda({\Bbb Z}_p,r)\stackrel{i}{\rightarrow} L_\lambda({\Bbb Z}_p)\stackrel{\pi}{\rightarrow}{\bf L}_\lambda^{[r]}({\Bbb Z}_p)
\rightarrow 0
$$
yields an exact sequence
$$
H^{i-1}(\Gamma,{\bf L}_\lambda^{[r]}({\Bbb Z}_p))\rightarrow H^i(\Gamma,L_\lambda({\Bbb Z}_p,r))\rightarrow H^i(\Gamma,L_\lambda({\Bbb Z}_p))
\stackrel{\pi^*}{\rightarrow} H^i(\Gamma,{\bf L}_\lambda^{[r]}({\Bbb Z}_p)).
$$
We denote by an upper index "TF" the maximal torsion-free quotient. Since $H^i(\Gamma,{\bf L}_\lambda^{[r]}({\Bbb Z}_p))$ is 
finite we further obtain an exact sequence
$$
0\rightarrow H^i(\Gamma,L_\lambda({\Bbb Z}_p,r))^{\rm TF}\stackrel{i^*}{\rightarrow} H^i(\Gamma,L_\lambda({\Bbb Z}_p))^{\rm TF}
\stackrel{\pi^*}{\rightarrow} Q\rightarrow 0,
$$
where $Q$ is a certain subquotient of $H^i(\Gamma,{\bf L}_\lambda^{[r]}({\Bbb Z}_p))$. 
%
%
%
%
%
%
%
%
%
%
%
%
%
%
%
%
%
Since $i^*$ is injective, we may identify $H^i(\Gamma,L_\lambda({\Bbb Z}_p,r))^{\rm TF}$ with its image under $i^*$. Hence, 
$$
\frac{H^i(\Gamma,L_\lambda({\Bbb Z}_p))^{\rm TF}}{H^i(\Gamma,L_\lambda({\Bbb Z}_p,r))^{\rm TF}}\quad\mbox{is a subquotient of}\;H^i(\Gamma,{\bf L}_\lambda^{[r]}({\Bbb Z}_p)).\leqno(3)
$$
We recall that $t$ equals the rank of $H^i(\Gamma,L_\lambda({\Bbb Z}_p))^{\rm TF}$ and we denote by 
$(p^{b_l})_l$, $b_1\ge b_2\ge \cdots\ge b_m>0$ the sequence of elementary divisors of 
$\frac{H^i(\Gamma,L_\lambda({\Bbb Z}_p))^{\rm TF}}{H^i(\Gamma,L_\lambda({\Bbb Z}_p,r))^{\rm TF}}$. Equations (2) and (3) yield $m\le n$ and
$$
b_1\le a_1,b_2\le a_2,\ldots,b_m\le a_m.\leqno(4)
$$
We set $b_l=0$ for $m<l\le t$ and $a_l=0$ for $n<l\le t$. 

\medskip

{\bf (1.5.2) } Using the results so far we can give a lower bound for ${\cal N}$. Since 
$H^i(\Gamma,L_{\tilde{\lambda}}({\Bbb Z}_p))\cong H^i(\Gamma,L_\lambda({\Bbb Z}_p))$ and 
$H^i(\Gamma,L_{\tilde{\lambda}}({\Bbb Z}_p,r))\cong H^i(\Gamma,L_\lambda({\Bbb Z}_p,r))$ as ${\Bbb Z}$-modules 
%
%
we obtain in particular, that
$$
\frac{H^i(\Gamma,L_{\tilde{\lambda}}({\Bbb Z}_p))^{\rm TF}}{H^i(\Gamma,L_{\tilde{\lambda}}({\Bbb Z}_p,r))^{\rm TF}}\cong \bigoplus_{l=1}^m \frac{\Bbb Z}{p^{b_l}{\Bbb Z}}
$$
is a finite abelian group and equations (2) and (4) imply that the $b_l$ are all smaller than or equal to $r$. Moreover, 
(5.4) Lemma 3.) in [M] and the definition of the action of the Hecke operator on cochains (cf. [M], section (5.2), equation (8)) 
imply that ${\Bbb T}(h)(Z^i(\Gamma,L_{\tilde{\lambda}}({\Bbb Z}_p,r)))\subseteq p^r Z^i(\Gamma,L_{\tilde{\lambda}}({\Bbb Z}_p))$, 
%
%
%
%
%
%
%
%
%
hence, we obtain 
$$
{\Bbb T}(h)(H^i(\Gamma,L_{\tilde{\lambda}}({\Bbb Z}_p,r))^{\rm TF})\subseteq p^r H^i(\Gamma,L_{\tilde{\lambda}}({\Bbb Z}_p))^{\rm TF}.
$$
%
%
%
Thus, we may apply Lemma 1 in [B] to the pair $L=H^i(\Gamma,L_{\tilde{\lambda}}({\Bbb Z}_p))^{\rm TF}$ and 
$K=H^i(\Gamma,L_{\tilde{\lambda}}({\Bbb Z}_p,r))^{\rm TF}$ and the operator $\xi={\Bbb T}(h)$. More 
precisely, we denote by $f_b:\,[0,t]\rightarrow {\Bbb R}_{\ge 0}$ the piecewise linear function attached to the sequence 
$(b_l)_{l=1,\ldots,t}$ as in [B]. We recall its definition: we set $C(j)=\sum_{l=1}^j (r-b_l)$, $0\le j\le t$, and $f_b$ is defined as the 
piecewise linear function joining the points $(j,C(j))$, $j=0,\ldots,t$. 
Lemma 1 in [B] then states that the Newton polygon ${\cal N}$ 
of ${\Bbb T}(h)$ acting on $H^i(\Gamma,L_{\tilde{\lambda}}({\Bbb Z}_p))^{\rm TF}$ is bounded from below by the graph of 
$f_b$. 
%
%

\medskip

{\bf (1.5.3) } We further estimate the function $f_b$. Equation (4) implies that $f_b$ lies above the piecewise linear function $f_a:\,[0,t]\rightarrow{\Bbb R}_{\ge 0}$ 
attached to the sequence $(a_l)_{l=1,\ldots,t}$ as defined in (1.5.2) (the ``$C(j)$'s'' {attached} to the sequence $(a_l)_l$ are smaller 
than those attached to the sequence $(b_l)_l$; note that also $a_l\le r$ for all $1\le l\le t$). Thus, we have
$$
{\cal N}\ge f_b\ge f_a.
$$
Using equation (2) it is not difficult to see that the function $f_a$ is the restriction to $[0,t]$ of the piecewise linear function 
which starts in $(0,0)$ and has slope $j$ for $g_i \sum_{h=0}^{j-1} N_h\le x\le g_i \sum_{h=0}^j N_h$, $j=0,\ldots,r-1$ and slope $r$ for 
$x\ge g_i \sum_{h=0}^{r-1} N_h$. 
%
%
By (1.1) we know that $N_h\le (h+1)^{s-1}$, $h\ge 0$; since also 
$g_i\le g$ we see that the function $f_a$ and, hence, ${\cal N}$, is bounded from below by the restriction to $[0,t]$ of the
piecewise linear function $f_r:\,{\Bbb R}_{\ge 0}\rightarrow{\Bbb R}$ which starts in $(0,0)$, has slope $j$ for 
$$
g \sum_{h=0}^{j-1} (h+1)^{s-1}\le x\le g \sum_{h=0}^j (h+1)^{s-1},\quad j=0,\ldots,r-1\leqno(5)
$$ 
and slope $r$ for $x\ge g \sum_{h=0}^{r-1} (h+1)^{s-1}$. Thus, we obtain
$$
{\cal N}\ge f_b\ge f_a\ge f_r.\leqno(6)
$$
This holds for all $r\in{\Bbb N}$ since $r$ was chosen arbitrarily. Taking into account that 
$$
\sum_{h=0}^{j-1} (h+1)^{s-1}=
\left\{
\begin{array}{ccc}
\frac{B_s(j+1)-B_s(0)}{s}&\mbox{if}&j\ge 1\\
0&\mbox{if}&j=0\\
\end{array}
\right.\leqno(7)
$$ 
it is immediate that $f_r$ and $f_\infty$ coincide on the interval $[0,g\frac{B_s(r+2)-B_s(0)}{s}]$. Since equation (6) in particular holds 
for arbitrarily large $r$ we finally obtain ${\cal N}\ge f_\infty$. 

\medskip

{\bf (1.5.4) } To finish the proof it remains to verify the alternative description of $f_\infty$. 
Equation (5) implies that $f_r$ connects the points 
$$
(0,0),\; P_j=(g \sum_{h=0}^j (h+1)^{s-1},g \sum_{h=0}^j h (h+1)^{s-1}),\quad j=0,1,\ldots,r-1
$$
(notice that the interval in equation (5) has length $(j+1)^{s-1}$ with corresponding slope $j$). Since $f_r$ and $f_\infty$ coincide on 
$[0,g\frac{B_s(r+2)-B_s(0)}{s}]$ for all $r\in{\Bbb N}$ this implies that $f_\infty$ is the piecewise linear function connecting the points 
$(0,0)$ and $P_j$, $j=0,1,2,\ldots$. Hence, the proof of the Theorem is complete.

\bigskip

{\bf (1.6) Corollary. }{\it For all $r\in{\Bbb N}$ the Newton polygon ${\cal N}$ lies above the restriction to $[0,t]$ of the 
piecewise linear function $f_\infty^\ast:\,{\Bbb R}_{\ge 0}\rightarrow{\Bbb R}_{\ge 0}$ which joins the points
$$
(0,0),\; Q_j=\left(g \frac{B_s(j+2)-B_s(0)}{s},g \frac{B_{s+1}(j+1)-B_{s+1}(0)}{s+1}\right), \quad j=0,1,2,\ldots
$$

}


\medskip

{\it Proof. } Again, using Bernoulli polynomials we obtain that the second coordinate of $P_j$ is bigger than or equal to 
$g \frac{B_{s+1}(j+1)-B_{s+1}(0)}{s+1}$, which is the second coordinate of $Q_j$. Hence, the claim follows. 

\bigskip

{\bf (1.7) Remark. } 1.) If $s\ge 2$ then $P_j$ lies strictly above $Q_j$ for all $j\ge 1$, hence, we obtain 
$f_\infty\ge f_\infty^\ast$. If $s=1$, i.e. in the rank $1$ case the points $P_j$ and $Q_j$ coincide, hence, $f_\infty=f_\infty^\ast$.

2.) The functions $f_\infty$ and $f_\infty^\ast$ do not depend on the weight ${\tilde{\lambda}}$.


\section{Explicit bounds for the dimension of the slope subspaces}

{\bf (2.1) } We apply the above Theorem to obtain bounds on the dimension $d=d({\tilde{\lambda}},i,\alpha)$ 
of the slope subspace $H^i(\Gamma,L_{\tilde{\lambda}}({\Bbb Q}_p))^{\le \alpha}$ and, hence, also on the dimension of 
$H^i(\Gamma,L_{\tilde{\lambda}}({\Bbb Q}_p))^\alpha$, $\alpha\in{\Bbb Q}_{\ge 0}$. To this end let $h:\,{\Bbb R}_{\ge 0}\rightarrow{\Bbb R}$ 
be any function such that $f_\infty^\ast(x)\ge h(x)$ for all $x\ge 0$.  
Since the function $x\mapsto (\alpha+\epsilon) x$, $\epsilon>0$, {\em strictly} lies above the Newton polygon ${\cal N}^{\le\alpha}={\cal N}|_{[0,d]}$ of 
${\Bbb T}(h)$ acting on $H^i(\Gamma,L_{\tilde{\lambda}}({\Bbb Q}_p))^{\le \alpha}$ 
(note that ${\cal N}^{\le\alpha}$ starts in $(0,0)$ and all its segments 
have slope $\le\alpha$)
we obtain from (1.6) Corollary 
$$
(\alpha+\epsilon) x > {\cal N}^{\le\alpha}(x)={\cal N}(x)\ge h(x)\quad\mbox{for all $0< x < d$}.
$$  
Thus, if $x\mapsto (\alpha+\epsilon) x$ and $h$ intersect in a point $(d(\epsilon),y)$ with $d(\epsilon)>0$, then necessarily $d\le d(\epsilon)$, 
hence, ${\rm inf}_{\epsilon>0} d(\epsilon)$ is an upper bound for $d=d({\tilde{\lambda}},i,\alpha)$.

\bigskip

{\bf (2.2) } We make an explicit choice of a function $h$ as in (2.1). To this end, let $P=x^s+a_{s-1}x^{s-1}+\cdots+a_0$ be a polynomial with coefficients 
$a_i\in{\Bbb R}$. For all $x\ge {\rm max}\{1,s|a_i|\}$ and all $i<s$ we have $|a_i|x^i\le |a_i| x^{s-1}\le \frac{1}{s}x^s$. Hence, we obtain
$$
P(x)\le 2x^s
$$
for all $x\ge {\rm max}\{1,s|a_{s-1}|,\ldots,s|a_0|\}$. Similarly, for all $x\ge {\rm max}\{1,2s|a_i|\}$ and all $i<s$ we have
$|a_i|x^i\le \frac{1}{2s}x^s$. Hence, we obtain
$$
P(x)\ge x^s-\sum_{i=0}^{s-1}|a_i| x^i\ge \frac{1}{2}x^s
$$
for all $x\ge {\rm max}\{1,2s|a_{s-1}|,\ldots,2s|a_0|\}$. Since $B_s$ is a polynomial of degree $s$ with leading coefficient equal to $1$ we
obtain as a special case of the above that there is $M=M(s)\in{\Bbb N}$ such that
$$
g\frac{B_s(x+2)-B_s(0)}{s}\le \frac{2g}{s} x^s \quad\mbox{and}\quad g\frac{B_{s+1}(x+1)-B_{s+1}(0)}{s+1}\ge \frac{g}{2(s+1)} x^{s+1}\leqno(1)
$$
for all $x\ge M$. The constant $M$ can be explicitly determined from the coefficients of the Bernoulli polynomials $B_s$ and $B_{s+1}$.

\medskip

We set $c=2^{-\frac{3s+1}{s}}(s/g)^{\frac{1}{s}}\in{\Bbb R}_{>0}$ and we define the 
function $h:\,[x_M,\infty)\rightarrow{\Bbb R}_{\ge 0}$, $x\mapsto c x^{\frac{s+1}{s}}$. Using equation (1) we obtain for all points $Q_j=(x_j,y_j)$ 
as in (1.6) Corollary with $j\ge M$:
$$
h(x_j)=h(g\frac{B_s(j+2)-B_s(0)}{s})\le  c (\frac{2g}{s} j^s)^{\frac{s+1}{s}}\le c  (\frac{2g}{s})^{\frac{s+1}{s}}\frac{2(s+1)}{g}y_j
\le y_j
$$
%
%
(notice that $2s\ge s+1$). Thus, $h$ lies below all points $Q_j=(x_j,y_j)$ with $j\ge M$ and since $h$ is convex and $f_\infty^\ast$ is the piecewise linear function 
connecting the points $Q_j$, the intermediate value Theorem implies that
$$
f_\infty^\ast(x)\ge h(x)\quad \mbox{for all $x\ge x_M=g \frac{B_s(M+2)-B_s(0)}{s}$}. \leqno(2)
$$
We extend $h$ to a function $h:\,{\Bbb R}_{\ge 0}\rightarrow{\Bbb R}$ by sending
$$
x\mapsto\left\{
\begin{array}{ccc}
f_\infty^\ast(x)&{\rm if}&x\le x_M\\   
h(x)&{\rm if}&x> x_M.\\
\end{array}
\right.
$$
and obtain $h(x)\le f^\ast_\infty(x)$ for all $x\ge 0$ by equation (2).

\bigskip

{\bf (2.3) Corollary. }{\it 1.) For all $\alpha\in{\Bbb Q}_{\ge 0}$, all dominant weights $\tilde{\lambda}\in X(\tilde{\bf T})$ and all $i$  we have
$$
{\rm dim}\,H^i(\Gamma,L_{\tilde{\lambda}}({\Bbb Q}_p))^{\le \alpha}\le m \alpha^s+ n,
$$
where $m=2^{3s+1}\frac{g}{s}$ and $n=g \frac{B_s(M+2)-B_s(0)}{s}$.

\medskip

2.) For all $\alpha\ge M(s)$, all dominant weights $\tilde{\lambda}\in X(\tilde{\bf T})$ and all $i$ we have
$$
{\rm dim}\,H^i(\Gamma,L_{\tilde{\lambda}}({\Bbb Q}_p))^{\le \alpha}\le m \alpha^s,
$$
where $m=2^{3s+1}\frac{g}{s}$.
}

\medskip

{\it Proof. } 1.) Since $\bullet$ $h(x)=f_\infty^\ast(x)=0$ in the interval $[0,g\frac{B_s(2)-B_s(0)}{s}]$ which has positive length by equation (7) in (1.5.3) 
$\bullet$ $h(x)$ grows to infinity faster than $x\mapsto(\alpha+\epsilon)x$ as 
$x\rightarrow\infty$ $\bullet$ $f^\ast_\infty(x_M)\ge h(x_M)$ (i.e. $h$ ``jumps downwards'' at $x_M$; cf. equation (2)) we see that the functions 
$x\mapsto(\alpha+\epsilon)x$ and $h(x)$, $x\ge 0$, always intersect in at least one point $S=(a,b)$ with $a>0$. By (2.1) $a$ then is a bound for the dimension 
$d=d(\tilde{\lambda},\alpha,i)$ of the slope $\alpha$ subspace. The functions $x\mapsto(\alpha+\epsilon)x$ and $x\mapsto c x^{\frac{s+1}{s}}$, $x\ge 0$, 
precisely intersect in the point $(d(\epsilon),y)$ with $d(\epsilon)=((\alpha+\epsilon)/c)^s$. We distinguish cases. If $(\alpha/c)^s\ge x_M$ then $d(\epsilon)\ge x_M$ 
and $(d(\epsilon),y)$ is a point of intersection of $x\mapsto(\alpha+\epsilon)x$ and $h$. Hence, we obtain $d\le d(\epsilon)$ for all $\epsilon>0$ which implies 
that $d\le (\alpha/c)^s=2^{3s+1}\frac{g}{s} \alpha^s$. Since $n\ge 0$ by equation (7) in (1.5.3) this implies the claim. If $(\alpha/c)^s<x_M$ then 
$x\mapsto(\alpha+\epsilon)x$ and $x\mapsto c x^{\frac{s+1}{s}}$ do not intersect in the range $x\ge x_M$. Hence, any intersection point $S=(a,b)$ of $h$ and 
$x\mapsto(\alpha+\epsilon)x$ satisfies $a\le x_M$ which implies $d\le x_M=g\frac{B_s(M+2)-B_s(0)}{s}$. Again, this yields the claim.

\medskip

2.) If $\alpha>M$ then $x\mapsto (\alpha+\epsilon) x$, $\epsilon>0$, has strictly bigger slope than any of the segments of $f_\infty$ in the intervall $[0,x_M]$ by (1.5) Theorem. Hence,  
$x\mapsto(\alpha+\epsilon) x$ strictly lies above $f_\infty^\ast$ 
and therefore above $h$ in the intervall $[0,x_M]$; the functions $x\mapsto(\alpha+\epsilon) x$ and $h$ 
therefore intersect in a point $S=(d(\epsilon),y)$ with $d(\epsilon)\ge x_M$. As in part 1.) we deduce that $d(\epsilon)=((\alpha+\epsilon)/c)^s$ which implies the claim. Hence,
the Corollary is proven.

\bigskip

{\bf Remark. } We did not make an attempt to optimize the coefficients $m,n$ appearing in (2.3) Corollary. This would mainly mean to replace the function $h$ by a 
function which is a better approximation from below to $f_\infty$ and, hence, to ${\cal N}$ or to improve on the bounds in equation (1). E.g. by allowing $n$ to become bigger 
one can achieve that $m$ comes arbitrarily close to $\frac{g}{s}$.

%
%
%
%
%
%
%
%
%
%

\bigskip

{\bf References}

\bigskip

[B-S] Borel, A., Serre, J.-P., Corners and arithmetic groups, Comment. Math. Helv. {\bf 48} (1974) 436 - 491

[B] Buzzard, K., Families of modular forms, J. de Th{\'e}orie de Nombres de Bordeaux {\bf 13} (2001) 43 - 52

[Br] Brown, Cohomology of groups, GTM {\bf 87}, Springer, 1982



[M] Mahnkopf, J., On truncation of irreducible representations of Chevalley groups, J. of Number Theory (2013)

[P] Pande, A., Local constancy of dimensions of Hecke eigenspaces of Automorphic forms, J. of Number Theory {\bf 129} (2009) 15 - 27

[R] Rotman, J., An introduction to the theory of groups, 4th ed., Springer 1995

[W] Wan, T., Dimension Variation of Classical and $p$-adic modular forms, Inv. Math. {\bf 133} (1998), 449 - 463.

\end{document}